\documentclass[12pt, a4paper, article, oneside]{memoir}

\usepackage[T1]{fontenc} 
\usepackage[latin1]{inputenc} 

\usepackage{framed}
\usepackage{hyperref}
\usepackage{doi}

\usepackage{amssymb}
\usepackage{amsmath, amsthm}

\usepackage{graphicx}
\usepackage{blindtext}
\usepackage{tabularray}
\UseTblrLibrary{booktabs}

\usepackage[listings,many]{tcolorbox} 
\usetikzlibrary{calc}
\usepackage{varwidth}
\usepackage{multirow}

\usepackage{tikz}
\usetikzlibrary{calc}
\usetikzlibrary[positioning, arrows]
\usetikzlibrary {intersections}
\usepackage{circuitikz}

\usepackage{soul,color,graphicx,colortbl} 
\usepackage{eurosym}
\usepackage{rotating,float}
\usepackage{hyperref} 
	\usepackage[figure]{hypcap}

\definecolor{cgray}{gray}{0.99} 
\definecolor{grayd}{rgb}{0.84,0.94,0.96} 
\definecolor{dgreen}{rgb}{0,0.5,0.25} 
\xdefinecolor{MyGreen}{cmyk}{0.92 ,0 ,0.87 ,0.09}

\definecolor{light-gray}{gray}{0.98}
\definecolor{background}{gray}{0.995}

\usepackage{soul, cancel, xcolor}

\definecolor{dred}{RGB}{191 0 64}

\usepackage{ulem}

\counterwithout{section}{chapter}
\setsecnumdepth{subsection}

\newtheorem{thm}{Theorem}

\newtheorem{example}{Example}%
\newtheorem{rem}{Remark}%

\newtheorem{prop}{Proposition}

\newtheorem{defn}{Definition}%

\newcommand{\uX}{\underline{x}}
\newcommand{\D}{\mathbf{D}_{\uX}}
\newcommand{\Dx}{\mathbf{D}_{(x_0,\uX)}}

\newcommand{\Jint}{\frak{A}\!\!\!\!\!\!\int}

\newcommand{\Dz}{\partial}
\newcommand{\Dzbar}{\overline{\partial_{q}}}

\begin{document}
	
\title{Conjugate $(1/q, q)$-harmonic Polynomials in $q$-Clifford Analysis}
\author{Swanhild Bernstein, Amedeo Altavilla, Martha Lina Zimmermann}
\date{}
\maketitle
\pagenumbering{arabic} 
\begin{abstract} We consider the problem of constructing a conjugate $(1/q, q)$-harmonic homogeneous polynomial $V_k$ of degree $k$ to a given $(1/q, q)$-harmonic homogeneous polynomial $U_k$ of degree $k.$ The conjugated harmonic polynomials $V_k$ and $U_k$ are associated to the $(1/q, q)$-mono\-genic polynomial $F = U_k + \overline{e}_0V_k. $
We investigate conjugate $(1/q, q)$-harmonic homogeneous polynomials in the setting of $q$-Clifford analysis. Starting from a given $(1/q, q)$-harmonic polynomial $U_k$ of degree $k$, we construct its conjugate counterpart $V_k$, such that the Clifford-valued polynomial $F = U_k + e_0 V_k$ is $(1/q, q)$-monogenic, i.e., a null solution of a generalized $q$-Dirac operator. Our construction relies on a combination of Jackson-type integration, Fischer decomposition, and the resolution of a $q$-Poisson equation. We further establish existence and uniqueness results, and provide explicit representations for conjugate pairs, particularly when $U_k$ is real-valued.

\end{abstract}
\section{Introduction}\label{Intro}

This work addresses the construction of conjugate harmonic polynomials within the framework of $q$-Clifford analysis, focusing on a specific $(1/q, q)$-Dirac operator. We begin in Section~1 with a motivation rooted in classical complex analysis, revisiting conjugate harmonic functions in the plane, introducing the foundations of $q$-calculus, and justifying our choice of operator. Section~2 provides the algebraic and analytic background, including a brief review of Clifford algebras and the formulation of $q$-deformed Dirac operators based on $1/q$- and $q$-partial derivatives. In Section~3, we develop the theory of homogeneous polynomials, introducing the $q$-Laplace operator, $q$-monogenic functions, and the associated Fischer decomposition. Finally, Section~4 presents our main results concerning conjugate $(1/q, q)$-harmonic polynomials and their relation to $(1/q, q)$-monogenic functions.

\subsection{Conjugate harmonic functions}
In the complex plane, a polynomial $p(x,y) = u(x,y) + i v(x,y)$ is holomorphic or analytic if $\partial_{\bar z} p=\frac{1}{2} (\partial_x + i \partial_y)p = 0$, which is equivalent to the Cauchy-Riemann system
\[
\left\{
\begin{array}{l}
\partial_x u = \partial_y v, \\
\partial_y u = - \partial_x v.
\end{array}
\right.
\]
With this holomorphic polynomial, we can associate the vector field $F(x,y) = (u(x,y), -v(x,y))^T$, which satisfies the Riesz system
\[
\left\{
\begin{array}{l}
\mathrm{div}\, F = 0, \\
\mathrm{curl}\, F = 0.
\end{array}
\right.
\]
If the vector field is conservative, there exists a real-valued harmonic potential $h(x,y)$ such that
\[
F = \mathrm{grad}\, h,
\]
which is related to the original holomorphic polynomial by 
\[
p = \partial_{z} h = \frac{1}{2}(\partial_x - i \partial_y)h.
\]

All these relations can be generalized to $n+1$ variables in Clifford analysis. Real-analytic polynomials and functions can be extended by the Cauchy-Kovalevskaya extension theorem to monogenic functions, i.e., null-solutions of the so-called Dirac operator. This operator is a first-order differential operator and serves as a higher-dimensional analogue of the $\partial_{\bar z}$-operator. It also plays a role analogous to the Cauchy-Riemann system.

If we consider a vector-valued monogenic function $F = \sum_{i=0}^n u_i e_i \sim (u_0, u_1, \ldots, u_n)^T \in \mathbb{R}^{n+1}$, then it satisfies $\mathrm{div}\, F = 0$ and $\mathrm{curl}\, F = 0$. On the other hand, the Riesz system can be extended to the Stein and Weiss system.

\subsection{$q$-Calculus}

Interestingly, $q$-calculus dates back to L.~Euler~\cite{Euler1748}. The $q$-derivative was introduced as an alternative to the ordinary derivative, even though it is, in fact, a difference operator. Important contributions to the field were made by F.H.~Jackson~\cite{Jackson1909}, H.~Exton~\cite{Exton1983}, and T.~Ernst, who has authored several papers and monographs on the subject~\cite{Ernst2003, Ernst2012}.

Quantum calculus--also referred to as $q$-calculus or $q$-deformation--is a framework in which one constructs $q$-analogues of classical operators, theorems, and identities. These analogues reduce to their classical counterparts as $q \to 1$. In this paper, we do not consider standard Clifford analysis but instead develop a $q$-Clifford analysis based on Jackson's formulation of $q$-calculus.

For a real number $u \in \mathbb{R}$ and a deformation parameter $q \in \mathbb{R},\, q\neq 1$, the $q$-deformation of $u$, sometimes called a {basic number}, is defined as
\begin{equation}
[u]_q = \frac{q^u - 1}{q - 1}.
\label{eq:qnumber}
\end{equation}
When $u \in \mathbb{N}$, equation~\eqref{eq:qnumber} simplifies to $[u]_q = 1 + q + q^2 + \cdots + q^{u-1}$. Taking the limit $q \to 1$, one recovers the classical identity $\lim_{q \to 1} [u]_q = u$.

Jackson calculus replaces the ordinary derivative with a difference operator:
\[
D_q f(x) = \frac{f(qx) - f(x)}{qx - x},
\]
which converges to the classical derivative in the limit $q \to 1$. In multiple dimensions, we introduce the $q$-partial derivatives
\[
\partial^q_{x_i} f(\uX) = \frac{\gamma_i f(\uX) - f(\uX)}{q x_i - x_i}, \quad \text{where } \gamma_i f(\uX) = f(x_1, \ldots, qx_i, x_{i+1}, \ldots, x_n).
\]

Functions such as analytic, harmonic, and monogenic functions play a central role in analysis. However, defining such functions in the context of $q$-calculus is not straightforward~\cite{Am}.

We will also require a notion of one-dimensional integration. Unfortunately, defining integration--or an \textit{antiderivative}--in $q$-calculus poses subtle challenges. As already discussed in~\cite{Kac2002}, we summarize the situation as follows.
Just like in classical calculus, an antiderivative is not unique. In the classical setting, the ambiguity is up to a constant, since the derivative of a function vanishes if and only if it is constant. In $q$-calculus, however, the condition $D_q \varphi(x) = 0$ is equivalent to $\varphi(qx) = \varphi(x)$, which does not necessarily imply that $\varphi$ is constant.

However, if we restrict our attention to formal power series, the situation becomes more manageable. If $\varphi(x) = \sum_{n=0}^\infty c_n x^n$ and $\varphi(qx) = \varphi(x)$, then this implies $q^n c_n = c_n$ for each $n$, which is possible only if $c_n = 0$ for all $n \geq 1$, i.e., $\varphi$ must be constant.
Therefore, if
\[
f(x) = \sum_{n=0}^\infty a_n x^n
\]
is a formal power series, then it admits a unique $q$-antiderivative (up to a constant term) in the space of formal power series, defined by
\[
\Jint f(x) \, \mathrm{d}_q x := \sum_{n=0}^\infty \frac{a_n x^{n+1}}{[n+1]_q} + C.
\]

Usually, one defines $\mathrm{d}_{q}x:=(q-1)x$. However, in this setting the symbol ``$\mathrm{d}_{q}x$'' only denotes the variable of \textit{antiderivation} with respect to the operator $D_{q}$.

Let $0 < q < 1$. Then, up to a constant, any function $f(x)$ has at most one $q$-antiderivative that is continuous at $x = 0$.
There also exists the notion of a Jackson integral \cite{Jackson1910} which is related to our notion of antiderivative used here. The Jackson integral for $q\in(0,1)$ is defined by
\begin{align*}
\int f(x) \,\mathrm{d}_q x = (1-q)\sum_{n=0}^\infty q^nxf(q^nx) + C.
\end{align*}
If $q\in(0,1)$ this series is a convergent geometric series. For polynomials it can be observed that this series converges to the expected function, we have
\begin{align*}
\int x^n \,\mathrm{d}_qx = \frac{1}{[n+1]_q}x^{n+1}+C.
\end{align*}
Because of these subtleties, it is often convenient to restrict attention to polynomials, and more specifically to homogeneous polynomials.

\subsection{$q$-Analytic Polynomials}

We aim to develop a $q$-calculus in which power functions and their inverses behave as holomorphic functions. This allows for the definition of Taylor and Laurent series in the $q$-setting~\cite{Pashaev2014}. In~\cite{Pashaev2014, Turner2016}, the authors introduce a notion of $q$-deformed complex analyticity for $0 < q < 1$. It is based on the operator 
\begin{equation} \label{Dzbar}
\Dzbar = \frac{1}{2} \left( \partial_x^q + i \partial_y^{1/q} \right),
\end{equation}
where $z = x + i y \in \mathbb{C}$, and the $q$-deformed partial derivatives are defined by
\begin{equation}
\partial_x^q f(x,y) = \frac{f(qx, y) - f(x, y)}{(q - 1)x}, \quad 
\partial_y^{1/q} f(x,y) = \frac{f(x, q^{-1}y) - f(x, y)}{(q^{-1} - 1)y}.
\end{equation}

As in the classical case, a function $f \colon D \subset \mathbb{C} \rightarrow \mathbb{C}$ is said to be $q$-analytic if $\Dzbar f = 0$. Given a real-analytic function, it is possible to extend it to a complex $q$-analytic function via an adapted version of the classical Cauchy-Kovalevskaya Theorem (see~\cite[Theorem 2]{Am}). In order to do so, we first introduce the following $q$-exponential function
\begin{defn}
For \( x \in \left(-\frac{1}{1 - q}, \frac{1}{1 - q}\right) \), the Jackson $q$-exponential function is defined as
\[
\exp_q(x) := \sum_{k=0}^{\infty} \frac{x^k}{[k]_q!}.
\]
\end{defn}

\begin{prop}[Cauchy-Kovalevskaya extension]
Let $p(x)$ be a real polynomial. Then
\begin{align*}
p(z) = p(x + i y) = \exp_{q^{-1}}(i y D_x^q)[f_0(x)] 
&= \sum_{k=0}^{\infty} \frac{1}{[k]_{1/q}!} (i y)^k (D_x^q)^k f_0(x) \\
&= \sum_{k=0}^{\infty} \frac{1}{[k]_q!} q^{(k-1)k/2} (i y)^k (D_x^q)^k f_0(x)
\end{align*}
is a $q$-analytic polynomial.
\end{prop}

\begin{example}
We show that standard real power functions extend naturally to the so-called $q$-binomials.
\begin{itemize}
\item For $f_0(x) = x$, we obtain
\[
F(x + i y) = \sum_{k=0}^{\infty} \frac{1}{[k]_q!} q^{(k-1)k/2} (i y)^k (D_x^q)^k x = x + i y.
\]

\item For $f_0(x) = x^{2}$, we obtain
\[
F(x + i y) = \sum_{k=0}^{\infty} \frac{1}{[k]_q!} q^{(k-1)k/2} (i y)^k (D_x^q)^k x^{2} = (x + i y)(x+iqy).
\]

\item For $f_0(x) = x^n$, we use the identity $(D_x^q)^k x^n = \frac{[n]_q!}{[n-k]_q!} x^{n-k}$ and find 
\begin{align*}
F(x + i y) &= \sum_{k=0}^{\infty} \frac{1}{[k]_q!} q^{(k-1)k/2} (i y)^k \frac{[n]_q!}{[n-k]_q!} x^{n-k} \\
&= \sum_{k=0}^{\infty} \binom{n}{k}_q q^{(k-1)k/2} (i y)^k x^{n-k} = (x + i y)(x+iqy)\cdots(x+iq^{n-1}y).
\end{align*}
\end{itemize}
\end{example}

Following this example the so-called complex $q$-binomials are defined by
\begin{equation}
z_q^k := (x + i y)_q^k = (x + i y)(x + i q y)(x + i q^2 y) \cdots (x + i q^{k-1} y),
\end{equation}
and serve as canonical examples of $q$-analytic functions~\cite{Pashaev2014}.

If we define $(x + i y)_q^{-k} := 1 / (x + i q^{-k} y)_q^k$, then we also have
\begin{equation}
\Dzbar (x + i y)_q^{-k} = 0, \quad 
\Dz (x + i y)_q^{-k} = [-k]_q (x + i y)_q^{-k-1}.
\end{equation}

\subsection{$q$-Analytic Functions are not Analytic Functions}

Since the operator $\Dzbar$ is not a {differential} operator, it is possible to find functions in its kernel that are not continuous.

\begin{example}[\cite{Am}] \label{Ex2}
Let $q \in (0,1) \cap \mathbb{Q}$, and define the function $\mathfrak{F} \colon \mathbb{C} \to \mathbb{C}$ by
\[
\mathfrak{F}(x + i y) = 
\begin{cases}
0, & \text{if } x \in \mathbb{Q}, \\
x + i y, & \text{if } x \in \mathbb{R} \setminus \mathbb{Q},
\end{cases}
\]
which satisfies $\Dzbar \mathfrak{F} = 0$. If $x \in \mathbb{Q}$, then $xq \in \mathbb{Q}$ as well, and thus $\mathfrak{F}(xq + i y) = \mathfrak{F}(x + i y) = 0$. If $x \in \mathbb{R} \setminus \mathbb{Q}$, then $xq \in \mathbb{R} \setminus \mathbb{Q}$, and the function behaves classically. Using similar constructions, one can generate infinitely many non-continuous $q$-analytic functions.
\end{example}

Therefore, $q$-analytic functions are best viewed as equivalence classes rather than functions in the classical sense. This makes it challenging to establish uniqueness principles. For this reason, we focus primarily on $q$-analytic functions that can be expressed as Taylor or Laurent series.

It is important to note that $q$-analytic functions are not analytic in the usual complex sense~\cite{Pashaev2014}. For instance, the $q$-binomials can be decomposed as
\[
z_q^k = z_{q^0} \cdot z_{q^1} \cdots z_{q^k},
\]
where
\begin{equation} \label{zqn}
z_{q^n} = x + i q^n y = \frac{1 + q^n}{2} z + \frac{1 - q^n}{2} \bar{z}.
\end{equation}
These expressions involve both $z$ and $\bar{z}$, placing $q$-analytic functions into the broader category of generalized analytic functions in the sense of Vekua~\cite{Vekua1962}.

\section{Preliminaries}\label{Pre}

\subsection{Classical Clifford Algebras}

We consider the Clifford algebra $\mathcal{C}\ell_{0,n+1}$ over $\mathbb{R}^{n+1}$, with identity element $1$ satisfying $1^2 = 1$, and generating elements $e_{0}, e_1, \ldots, e_n$ fulfilling the multiplication rules
\begin{equation*}
	e_i e_j + e_j e_i = -2 \delta_{ij}
\end{equation*}
for $i, j = 0, 1, \ldots, n$. The elements $e_i$ generate a basis $\mathcal{B}_0$ consisting of $2^{n+1}$ elements:
\[
\mathcal{B}_0 = \left(1,e_0, e_1, e_2, \ldots, e_0 e_1, \ldots, e_0 e_1 e_2, \ldots \right).
\]

Let $M := \{0, \ldots, n\}$ and define $A := \{(h_1, \ldots, h_r) \in \mathcal{P}M : 0 \leq h_1 \leq \cdots \leq h_r \leq n\}$. An arbitrary element $\lambda \in \mathcal{C}\ell_{0,n}$ can be written as
\begin{equation*}
	\lambda = \sum_A \lambda_A e_A, \quad \lambda_A \in \mathbb{R},
\end{equation*}
with $e_A = e_{h_1} \cdots e_{h_r}$. The length (or norm) of a Clifford number is defined as 
\begin{equation*}
	|\lambda|_0 := 2^{(n+1)/2} \left( \sum_A |\lambda_A|^2 \right)^{1/2}.
\end{equation*}

A Clifford conjugation in $\mathcal{C}\ell_{0,n+1}$ is defined by $\overline{e}_i = -e_i$ for $i = 0, 1, \ldots, n$, extended linearly so that $\overline{e_A e_B} = \overline{e}_B \, \overline{e}_A$ and $\overline{e}_{\emptyset} = e_{\emptyset}$.

The Clifford algebra $\mathcal{C}\ell_{0,n+1}$ admits the decomposition
\begin{equation} \label{split}
\mathcal{C}\ell_{0,n+1} = \mathcal{C}\ell_{0,n} + \overline{e}_0 \, \mathcal{C}\ell_{0,n},
\end{equation}
where $\mathcal{C}\ell_{0,n}$ is the subalgebra generated by $\{e_1, \ldots, e_n\}$.

\subsection{Clifford Analysis}

We consider $\mathcal{C}\ell_{0,n+1}$-valued polynomials $P$ on $\mathbb{R}^{n+1}$:
\begin{equation*}
	P(x_0, \underline{x}) = \sum_A P_A(x_0, \underline{x}) e_A = \sum_A P_A(x_0, x_1, \ldots, x_n) e_A, \quad \text{with } P_A \colon \mathbb{R}^{n+1} \rightarrow \mathbb{R},
\end{equation*}
with $A\subset\mathcal{B}_0$ and $\mathcal{C}\ell_{0,n}$-valued polynomials $p$ on $\mathbb{R}^n$: 
\begin{equation*}
	p(\underline{x}) = \sum_C p_C(\underline{x}) e_C = \sum_C p_C(x_1, \ldots, x_n) e_C, \quad \text{with } p_C \colon \mathbb{R}^n \rightarrow \mathbb{R}.
\end{equation*}
Here, $C\subset \mathcal{B}$ with 
\begin{align*}
\mathcal{B} = \left(1, e_1, e_2, \ldots,  e_1e_2, \ldots,  e_1 e_2e_3, \ldots \right)
\end{align*}
and $C := \{(h_1, \ldots, h_r) \in \mathcal{P}K : 0 \leq h_1 \leq \cdots \leq h_r \leq n\}$ for $K := \{1,\ldots,n\}$.
A vector in $\mathbb{R}^n$ can be identified with a Clifford-valued variable:
\begin{equation*}
	\underline{x} = \sum_{i=1}^n x_i e_i.
\end{equation*}
From the multiplication rules, we deduce that $\underline{x}^2$ is scalar-valued:
\[
\underline{x}^2 = -\sum_{i=1}^n x_i^2 = -|\underline{x}|^2.
\]

For more on Clifford analysis, see the monographs~\cite{BDS1982, DSS1992, GM1991}.
An initial approach to $q$-calculus in Clifford analysis was proposed in~\cite{CoSo2010, CoSo2011}. Our method, presented in~\cite{BeZiSchn2022}, follows a different direction. A more general framework for $q$-Clifford analysis was developed in~\cite{GonzalezCervantes2024}, where a Borel-Pompeiu formula was also established~\cite{GonzalezCervantes2024a}.

We now introduce $q$-Clifford analysis by defining two Dirac-type operators.

First, the $q$-Dirac operator $\D^q$ on $\mathbb{R}^n$ is defined as
\begin{equation*}
	\D^q = \sum_{i=1}^n e_i \, \partial^q_{x_i}.
\end{equation*}
This operator satisfies $(\D^q)^2 = -\Delta^q_{\underline{x}}$, where $\Delta^q_{\underline{x}}$ is the $q$-Laplace operator on $\mathbb{R}^n$.

Second, we define the Dirac operator on $\mathbb{R}^{n+1}$, analogous to the complex operator \eqref{Dzbar}, as
\[
\Dx^{1/q, q} = e_0 \, \partial^{1/q}_{x_0} + \D^q.
\]
This leads to an associated Laplace-type operator given by
\begin{equation}
  \Dx^{1/q, q} \Dx^{1/q, q} = - (\partial^{1/q}_{x_0})^2 + (\D^q)^2 = -( (\partial^{1/q}_{x_0})^2 + \Delta^q_{\uX}).
\end{equation}

\section{Homogeneous harmonic polynomials}

\subsection{$q$-Harmonics and $q$-Monogenics}

We start with the definition of a homogeneous polynomial. A polynomial $P(\uX)$ is called homogeneous of degree $k$ if, for every $\uX \in \mathbb{R}^m$ with $\uX \neq 0$, it satisfies
\begin{equation*}
P(\uX) = \vert \uX \vert^k P\left(\frac{\uX}{\vert \uX \vert}\right).
\end{equation*}

The $q$-harmonic functions are the null-solutions of the $q$-Laplace operator $\Delta^q_{\uX} = \sum_{i=1}^n (\partial^q_{x_i})^2$. The space of $k$-homogeneous polynomials that are also $q$-harmonic is denoted by
\begin{equation*}
\mathcal{H}_k = \{ P(\uX) \in \mathcal{P}_k : \Delta^q_{\uX} P(\uX) = 0 \}.
\end{equation*}

We now consider the Clifford-valued case. A $\mathcal{C}\ell_{0,n}$-valued homogeneous polynomial can be expressed as
\begin{align*}
P(\uX) = \sum_C c_C \uX^C, \quad c_C \in \mathcal{C}\ell_{0,n} \quad \Leftrightarrow \quad P(\uX) = \sum_C e_C p_C(\uX),
\end{align*}
where each $p_C(\uX)$ is a real-valued $k$-homogeneous polynomial and $|C| = k$.

As previously noted, monogenic functions are a central object of study in Clifford analysis. The space of monogenic homogeneous polynomials of degree $k$ is denoted by
\begin{equation*}
\mathcal{M}_k = \{ M(\uX) \in \mathcal{P}_k : \D^q M = 0 \}.
\end{equation*}

\subsection{Fischer Decomposition on $\mathbb{R}^n$}

We use multi-index notation $\alpha = (\alpha_0, \alpha_1,\ldots ,\alpha_n) \in \mathbb{N}^{n+1}$, with:
\begin{itemize}
\item $\uX^{\alpha} = x_1^{\alpha_1} \cdots x_m^{\alpha_m}$ and $x^{\alpha} = x_0^{\alpha_0} x_1^{\alpha_1} \cdots x_m^{\alpha_m}$,
\item $\alpha! = \alpha_0! \alpha_1! \cdots \alpha_n!$,
\item $|\alpha| = \sum_{i=0}^n \alpha_i$,
\item $(\D^q)^{\alpha} = (\partial_{x_0}^q)^{\alpha_0} (\partial_{x_1}^q)^{\alpha_1} \cdots (\partial_{x_m}^q)^{\alpha_m}$.
\end{itemize}

We consider homogeneous polynomials $p(\uX)$ in $\mathbb{R}^n$. The standard basis for $k$-homogeneous Clifford-valued polynomials is $\mathcal{P}_k = \{ \uX^{\alpha} : |\alpha| = k \}$.

We now define an inner product on the complex vector space $\mathcal{P}_k$. For $R_1, R_2 \in \mathcal{P}_k$, where each $R_i(\uX) = \sum_{|\alpha| = k} x^{\alpha} a_{\alpha}^i$ with $a_{\alpha}^i \in \mathcal{C}\ell_{0,m}$ for $i = 1, 2$, we define the Fischer inner product as
\begin{equation*}
\langle R_1, R_2 \rangle_{k,q} = \sum_{|\alpha| = k} [\alpha]_q! \, (\overline{a_{\alpha}^1} a_{\alpha}^2)_0.
\end{equation*}

Propositions~\ref{thm:fischerip} through~\ref{thm:fischerdecomp} are proven in~\cite{BeZiSchn2022} for monogenic functions. The case of harmonic functions follows analogously.

\begin{prop}
For $R_1, R_2 \in \mathcal{P}_k$, the Fischer inner product satisfies
\begin{equation*}
\langle R_1, R_2 \rangle_{k,q} = (\overline{R}_1(\D^q) R_2)_0,
\end{equation*}
where $\overline{R}_1(\D^q)$ denotes the operator obtained by replacing $x_j$ in $R_1$ with the $q$-partial derivative $\partial_{x_j}^q$.
\label{thm:fischerip}
\end{prop}

\begin{prop}
For all $Q \in \mathcal{P}_k$, $P \in \mathcal{P}_{k+1}$, and $R \in \mathcal{P}_{k+2}$, we have
\begin{equation}
\langle \uX Q, P \rangle_{k+1,q} = - \langle Q, \D^q P \rangle_{k,q}, \quad 
\langle -|\uX|^2 Q, R \rangle_{k+2,q} = \langle Q, \Delta^q_{\uX} R \rangle_{k,q}.
\end{equation}
\label{thm:fischerip2}
\end{prop}

\begin{prop}
For $k \in \mathbb{N}$, the following decomposition holds:
\begin{equation}
\mathcal{P}_k = \mathcal{M}^q_k \oplus \uX \mathcal{P}_{k-1} = \mathcal{H}^q_k \oplus |\uX|^2 \mathcal{P}_{k-2}.
\end{equation}
Moreover, the subspaces $\mathcal{M}^q_k$ and $\uX \mathcal{P}_{k-1}$, as well as $\mathcal{H}^q_k$ and $|\uX|^2 \mathcal{P}_{k-2}$, are orthogonal with respect to the Fischer inner product.
\label{thm:fischerdecomp}
\end{prop}

We therefore obtain the Fischer decomposition of the space of homogeneous polynomials $\mathcal{P}_k$:
\begin{equation}
\mathcal{P}_k = \sum_{s=0}^{k} \uX^s \mathcal{M}^q_{k-s}, \qquad 
\mathcal{P}_k = \sum_{s=0}^{\left\lfloor \frac{k}{2} \right\rfloor} |\uX|^{2s} \mathcal{H}^q_{k-2s}.
\end{equation}

Furthermore, the space of $q$-harmonic polynomials of degree $k$, denoted by $H_k$, can be decomposed into $q$-monogenic components:
\begin{equation*}
H_k = M_k + \uX M_{k-1}.
\end{equation*}

\section{$(1/q,q)$-Conjugate Harmonic Polynomials in $q$-Clifford Analysis}

In this section, we generalize the results of~\cite{BD2003, BD2003C}, where conjugate harmonic functions were considered in the context of Clifford analysis.

The basis for conjugate harmonic polynomials is the decomposition~\eqref{split} of the Clifford algebra $\mathcal{C}\ell_{0,n+1}$. Accordingly, we split a $\mathcal{C}\ell_{0,n+1}$-valued polynomial $F$ into $\mathcal{C}\ell_{0,n}$-valued components $U$ and $V$:
\[
F(x_0, \uX) = U(x_0, \uX) + \overline{e}_0 V(x_0, \uX).
\]

\begin{defn}
A $\mathcal{C}\ell_{0,n+1}$-valued polynomial $P = U + \overline{e}_0 V$ is called $(1/q, q)$-monogenic in $\mathbb{R}^{n+1}$ if and only if the $\mathcal{C}\ell_{0,n}$-valued polynomials $U$ and $V$ satisfy the system:
\begin{align}
\partial^{1/q}_{x_0} U + \D^q V &= 0, \label{qconharm1} \\
\D^q U + \partial^{1/q}_{x_0} V &= 0. \label{qconharm2}
\end{align}
\end{defn}

It follows immediately that both $U$ and $V$ are $(1/q, q)$-harmonic polynomials.

\begin{defn}
A pair $(U, V)$ of $\mathcal{C}\ell_{0,n}$-valued $(1/q, q)$-harmonic polynomials in $\mathbb{R}^{n+1}$ is called a pair of conjugate $(1/q, q)$-harmonic polynomials if $P = U + \overline{e}_0 V$ is $(1/q, q)$-monogenic in $\mathbb{R}^{n+1}$.
\end{defn}
As mentioned in Section 1.3, the Jackson integral only converges for $q\in(0,1)$. Now we will integrate with respect to $\frac {1}{q}$; therefore, from now on we consider the case $q>1$.
\begin{rem}
We will need to integrate a homogeneous polynomial $p(x_0, \uX) = x_0^{\alpha_0} x_1^{\alpha_1} \cdots x_n^{\alpha_n}$ with respect to $x_0$. The antiderivative with respect to $\partial^{1/q}_{x_0}$ is denoted by $\mathfrak{A} \!\!\!\!\!\int$, and satisfies:
\[
\Jint_0^{x_0} p(t, \uX)\, d_{\frac{1}{q}}t = \frac{1}{[\alpha_0+1]_{1/q}} x_0^{\alpha_0 + 1} x_1^{\alpha_1} \cdots x_n^{\alpha_n}.
\]
\end{rem}

Let $U_k(x_0, \uX)$ be a $\mathcal{C}\ell_{0,n+1}$-valued $(1/q, q)$-harmonic homogeneous polynomial of degree $k$ in $\mathbb{R}^{n+1}$. We aim to construct a $\mathcal{C}\ell_{0,n}$-valued homogeneous polynomial $V_k(x_0, \uX)$ of degree $k$ that is conjugate $(1/q, q)$-harmonic to $U_k$.

Since $V_k$ is to be conjugate to $U_k$, the pair must satisfy the system~\eqref{qconharm1}--\eqref{qconharm2}. We begin by integrating~\eqref{qconharm2} and observe that in~\eqref{qconharm1}, $V_k$ is only determined up to a $q$-harmonic function. We make the ansatz:
\begin{equation}\label{ansatz}
V_k(x_0, \uX) = - \Jint_0^{x_0} \D^q U_k(t, \uX)\, d_{\frac{1}{q}}t + W(\uX),
\end{equation}
where $W(\uX)$ is a $\mathcal{C}\ell_{0,n}$-valued polynomial on $\mathbb{R}^n$.

Recalling that $U_{k}$ is $(1/q,q)$-harmonic, we compute:
\begin{align*}
\D^q V_k(x_0, \uX) &= - \Jint_0^{x_0} (\D^q)^2 U_k(t, \uX)\, d_{\frac{1}{q}}t + \D^q W(\uX) \\
&= - \Jint_0^{x_0} (\partial^{1/q}_{x_0})^2 U_k(t, \uX)\, d_{\frac{1}{q}}t + \D^q W(\uX) \\
&= -\partial^{1/q}_{x_0} U_k(x_0, \uX) + \partial^{1/q}_{x_0} U_k(0, \uX) + \D^q W(\uX).
\end{align*}
Thus, for~\eqref{qconharm1} to hold, we require:
\[
\D^q W(\uX) = - \partial^{1/q}_{x_0} U_k(0, \uX).
\]

Since $U_k(x_0, \uX)$ is a homogeneous polynomial of degree $k$ in $(x_0, \uX)$, it can be written as
\[
U_k(x_0, \uX) = \sum_A \left( \sum_{|\alpha| = k} a_{A,\alpha} x^{\alpha} \right) e_A,
\]
where $x^\alpha = x_0^{\alpha_0} x_1^{\alpha_1} \cdots x_n^{\alpha_n}$. Therefore, the derivative $\partial^{1/q}_{x_0} U_k(0, \uX)$ is either zero or a homogeneous polynomial of degree $k-1$ in $\uX$.

If $\partial^{1/q}_{x_0} U_k(0, \uX) = 0$, we choose $W(\uX) = 0$. Otherwise, Fischer decomposition guarantees the existence of a $\mathcal{C}\ell_{0,n}$-valued homogeneous polynomial $h_{k-1}(\uX)$ of degree $k-1$ such that
\[
h_{k+1}(\uX) = -|\uX|^2 h_{k-1}(\uX),
\]
where $h_{k+1}$ is the unique homogeneous polynomial of degree $k+1$ satisfying
\[
\Delta^q_{\uX} h_{k+1}(\uX) = h_{k-1}(\uX).
\]
Combining this with the previous identity, we see that there always exists a $\mathcal{C}\ell_{0,n}$-valued polynomial $h_{k+1}(\uX)$ such that
\begin{equation} \label{401}
\Delta^q_{\uX} h_{k+1}(\uX) = \partial^{1/q}_{x_0} U_k(0, \uX).
\end{equation}

If we take  
\[
W(\uX) = \D^q h_{k+1}(\uX),
\]
then \( W(\uX) \) is a $\mathcal{C}\ell_{0,n}$-valued homogeneous polynomial of degree $k$ in $\mathbb{R}^n$ and satisfies
\[
\D^q W(\uX) = (\D^q)^2 h_{k+1}(\uX) = -\Delta^q_{\uX} h_{k+1}(\uX) = - \partial^{1/q}_{x_0} U_k(0, \uX).
\]

Hence, for a given $\mathcal{C}\ell_{0,n+1}$-valued $(1/q, q)$-harmonic polynomial $U_k(x_0, \uX)$ of degree $k$ in $\mathbb{R}^{n+1}$, we obtain:
\begin{align}
V_k(x_0, \uX) &= - \Jint_0^{x_0} \D^q U_k(t, \uX)\, d_{\frac{1}{q}}t + \D^q h_{k+1}(\uX) \nonumber \\
&= -\D^q \left( \Jint_0^{x_0} U_k(t, \uX)\, d_{\frac{1}{q}}t - h_{k+1}(\uX) \right), \label{vk}
\end{align}
where \( h_{k+1} \) satisfies the $q$-Poisson equation:
\[
\Delta^q_{\uX} h_{k+1}(\uX) = \partial^{1/q}_{x_0} U_k(0, \uX).
\]

This completes the construction of the conjugate $(1/q, q)$-harmonic polynomial \( V_k \) associated with \( U_k \), and we summarize this result in the following theorem.

\begin{thm}
Let $U_k(x_0, \uX)$ be a $\mathcal{C}\ell_{0,n+1}$-valued $(1/q, q)$-harmonic function in $\mathbb{R}^{n+1}$, and let $h_{k+1}(\uX)$ satisfy the $q$-Poisson equation~\eqref{401}. Define
\begin{equation} \label{H}
H_{k+1}(x_0, \uX) = \Jint_0^{x_0} U_k(t, \uX)\, d_{\frac{1}{q}}t - h_{k+1}(\uX).
\end{equation}
Then:
\begin{itemize}
\item[(i)] $H_{k+1}$ is $\mathcal{C}\ell_{0,n}$-valued and $(1/q, q)$-harmonic in $\mathbb{R}^{n+1}$.
\item[(ii)] $V_k(x_0, \uX) = -\D^q H_{k+1}(x_0, \uX)$ is $\mathcal{C}\ell_{0,n}$-valued and conjugate $(1/q, q)$-harmonic to $U_k$.
\item[(iii)] $F = U_k + \overline{e}_0 V_k = \Dx^{1/q,q} \overline{e}_0 H_{k+1}$ is $(1/q, q)$-monogenic in $\mathbb{R}^{n+1}$.
\end{itemize}
\end{thm}

Hence,
\[
H_{k+1}(x_0, \uX) = \Jint_0^{x_0} U_k(t, \uX)\, d_{\frac{1}{q}}t - h_{k+1}(\uX)
\]
is a $\mathcal{C}\ell_{0,n}$-valued polynomial in $\mathbb{R}^{n+1}$ such that
\begin{align*}
F(x_0, \uX) &= U_k(x_0, \uX) + \overline{e}_0 V_k(x_0, \uX) \\
&= \partial^{1/q}_{x_0} H_{k+1}(x_0, \uX) + \overline{e}_0 (-\D^q H_{k+1}(x_0, \uX)) \\
&= (e_0 \partial^{1/q}_{x_0} + \D^q) \overline{e}_0 H_{k+1}(x_0, \uX).
\end{align*}

Since \( \Dx^{1/q, q} F = 0 \) in \( \mathbb{R}^{n+1} \), it follows that \( \overline{e}_0 H_{k+1}(x_0, \uX) \) is \( (1/q, q) \)-harmonic in \( \mathbb{R}^{n+1} \). \hfill $\bullet$

This leads to the following structural result. Let \( M^+(k; \mathcal{C}\ell_{0,n+1}) \) denote the space of \( \mathcal{C}\ell_{0,n} \)-valued homogeneous \( (1/q, q) \)-harmonic polynomials of degree \( k \) in \( \mathbb{R}^{n+1} \). Then:

\begin{thm}
Let \( P_k \in M^+(k; \mathcal{C}\ell_{0,n+1}) \). Then there exists a \( \mathcal{C}\ell_{0,n} \)-valued homogeneous \( (1/q, q) \)-harmonic polynomial \( H_{k+1} \) of degree \( k+1 \) in \( \mathbb{R}^{n+1} \) such that
\[
P_k = \overline{e}_0 (e_0 \partial^q_{x_0} - \D^q) H_{k+1}.
\]
\end{thm}

A special case arises when \( u_k(x_0, \uX) \) is a \textbf{real-valued} homogeneous \( (1/q, q) \)-harmonic polynomial. From equations~\eqref{401} and~\eqref{H}, it follows that both \( h_{k+1}(\uX) \) and \( H_{k+1}(x_0, \uX) \) are also real-valued. Thus, by~\eqref{vk}, the corresponding conjugate \( v_k(x_0, \uX) \) is \( \mathcal{C}\ell_{0,n} \)-valued. This leads to

\begin{thm} \label{T3}
Given a real-valued homogeneous \( (1/q, q) \)-harmonic polynomial \( u_k(x_0, \uX) \) of degree \( k \) in \( \mathbb{R}^{n+1} \), there exists a unique 
\( \mathcal{C}\ell_{0,n} \)-valued homogeneous \( (1/q, q) \)-harmonic polynomial \( v_k(x_0, \uX) \) of degree \( k \), conjugate to \( u_k(x_0, \uX) \), of the form
\[
v_k(x_0, \uX) = -\D^q H_{k+1}(x_0, \uX), \quad \text{where} \quad H_{k+1}(x_0, \uX) = \Jint_0^{x_0} u_k(t, \uX)\, d_{\frac{1}{q}}t - h_{k+1}(\uX).
\]
\end{thm}

\begin{rem}
Using the Leibniz rule in $q$-calculus:
\[
\partial^q_{x_i} (p_1(\uX) p_2(\uX)) = \partial^q_{x_i}(p_1(\uX)) \, \gamma_i(p_2(\uX)) + p_1(\uX) \, \partial^q_{x_i}(p_2(\uX)),
\]
we find:
\[
e_i \partial^q_{x_i} \left( |\uX|^2 h_{k-1}(\uX) \right) = [2]_q x_i e_i \gamma_i(h_{k-1}(\uX)) + |\uX|^2 e_i \partial^q_{x_i} h_{k-1}(\uX).
\]

Therefore, the conjugate polynomial \( v_k \) in Theorem~\ref{T3} can be rewritten as
\[
v_k(x_0, \uX) = v^{(1)}_k(x_0, \uX) + \sum_{i=1}^n x_i e_i \gamma_i(w^{(1)}_{k-1}(\uX)) + |\uX|^2 w^{(2)}_{k-2}(\uX),
\]
where:
\begin{itemize}
\item[(i)] \( \displaystyle v^{(1)}_k(x_0, \uX) = \Jint_0^{x_0} \D^q u_k(t, \uX)\, d_{\frac{1}{q}}t \) is a homogeneous polynomial of degree \( k \) in \( \mathbb{R}^{n+1} \),
\item[(ii)] \( w^{(1)}_{k-1}(\uX) = [2]_q h_{k-1}(\uX) \) is a homogeneous polynomial of degree \( k-1 \) in \( \mathbb{R}^n \),
\item[(iii)] \( w^{(2)}_{k-2}(\uX) = -\D^q h_{k-1}(\uX) \) is a homogeneous polynomial of degree \( k-2 \) in \( \mathbb{R}^n \).
\end{itemize}
\end{rem}

\begin{example}
Let us compute an explicit example. Set \( q = \frac{4}{3} \) (so that \( \frac{1}{q} = \frac{3}{4} \)) and consider the polynomial
\[
U_3(x_0, x_1, x_2) = x_0^3 - x_0 x_1^2 - \frac{47}{64} x_0 x_2^2.
\]
We compute the \((1/q, q)\)-Laplace operator:
\[
(\partial_{x_0}^{\frac{1}{q}})^2 U_3 = \partial_{x_0}^{\frac{1}{q}}\left([3]_{\frac{1}{q}} x_0^2 - x_1^2 - \frac{47}{64} x_2^2\right) = [3]_{\frac{1}{q}} [2]_{\frac{1}{q}} x_0,
\]
and
\[
(\partial_{x_1}^{q})^2 U_3 = -[2]_{q} x_0, \qquad (\partial_{x_2}^{q})^2 U_3 = -[2]_{q} \cdot \frac{47}{64} x_0.
\]
Therefore, the total \((1/q, q)\)-Laplacian vanishes, and \( U_3 \) is \((1/q, q)\)-harmonic.

From the ansatz~\eqref{ansatz}, we obtain
\[
V_3(x_0, \uX) = -[2]_{q} x_0^2 \left(1 + \frac{47}{64} \right) + W(\uX),
\]
where the correction term \( W(\uX) \) solves
\[
\D^q W(\uX) = -\partial_{x_0}^{\frac{1}{q}} U_3(0, \uX) = -x_1^2 - \frac{47}{64} x_2^2.
\]
A suitable solution is given by
\[
h_4(\uX) = -\frac{1}{[4]_q [3]_q} \left( x_1^4 + \frac{47}{64} x_2^4 \right),
\]
which satisfies
\[
\Delta^q_{\uX} h_4(\uX) = \partial_{x_0}^{\frac{1}{q}} U_3(0, \uX).
\]
Then the corresponding correction term is
\[
W(\uX) = \D^q h_4(\uX) = -e_1 \frac{x_1^3}{[3]_q} - e_2 \frac{47}{64} \cdot \frac{x_2^3}{[3]_q},
\]
which completes the construction of \( V_3 \).

\end{example}

\section*{Acknowledgement}
Amedeo Altavilla was partially supported by PRIN 2022MWPMAB - ``Interactions between Geometric Structures and Function Theories'' and by GNSAGA of INdAM.

\noindent Swanhild Bernstein: \url{swanhild.bernstein@math.tu-freiberg.de} \\
Amedeo Altavilla: \url{amedeo.altavilla@uniba.it} \\
Martha Lina Zimmermann: \url{Martha-Lina.Zimmermann@math.tu-freiberg.de}

\end{document}